\documentclass[conference]{IEEEtran}

\usepackage{pgfplots,pgfplotstable,booktabs}

\pgfplotsset{compat=1.5}

\usepackage{amsmath, bbold}

\usepackage{etoolbox}

\usepackage{graphicx, epstopdf}
\usepackage{caption}
\usepackage{subfigure}

\DeclareGraphicsExtensions{.eps}

\apptocmd{\sloppy}{\hbadness 10000\relax}{}{}

\usepackage{cite}

\usepackage{epsf}

\begin{document}

\title{Simple Certificate of Solvability of Power Flow Equations for Distribution Systems}

\author{
\IEEEauthorblockN{Suhyoun Yu, Hung D. Nguyen, and Konstantin S. Turitsyn}
\IEEEauthorblockA{Department of Mechanical Engineering\\
Massachusetts Institute of Technology\\
Cambridge, MA 02139\\
Email: syu2@mit.edu, hunghtd@mit.edu, and turitsyn@mit.edu}}

% make the title area
\maketitle

\begin{abstract}
%\boldmath
Power flow solvable boundary plays an important role in contingency analysis, security assessment, and planning processes. However, to construct the real solvable boundary in multidimensional parameter space is burdensome and time consuming. In this paper, we develop a new technique to approximate the solvable boundary of distribution systems based on Banach fixed point theorem. Not only the new technique is fast and non-iterative, but also the approximated boundary is more valuable to system operators in the sense that it is closer to the feasible region. Moreover, a simple solvable criterion is also introduced that can serve as a security constraint in various planning and operational problems.

\end{abstract}


\section{Introduction}

% {\color{blue}
% % @Hung (from essie) Line 86: are there supposed to be some equations? No. I think Sharelarex made the gap itself so that section II begins at the second page.\\
% @Prof. Turitsyn: From \eqref{eq:solcri}, if we introduce a new variable: $\sigma = \sum_k |s_k/s_{k,max}|$. Could $1-\sigma$ be used as a security indicator?\\
% % @Essie: Does $s_{k,max}$ in \eqref{eq:solcri} depend on the slack bus voltage $V_0$? or you assume that $V_0 = 1$.
% }

The penetration of renewable energy sources is expected to increase substantially in the near future. Besides numerous advantages of such more environmentally friendly resources, additional risks may make the distribution system more vulnerable to instability \cite{KundurDSA}. Voltage security is one of main existing issues that we need to address carefully and thoroughly in order to operate the systems in a reliable and secure fashion. Hence, the system operators also need to know the system limits which require faster and more reliable security assessment tools. In fact, static voltage stability criteria is widely used by system operators \cite{xie2007novel,LeXieQSVS}; moreover, it has been argued that static analysis is preferred over dynamic approach \cite{morison1993voltage}. Most of static analysis techniques rely on the solution of quasi-steady state analysis or power flow problem. The concept of quasi-steady state is based on the assumption of slow dynamics so that the power system can be considered in steady state.

% Power flow in a network is determined by the voltage at each bus of the network and the impedance of the lines between buses. The load flow problem consists of finding the set of voltages: magnitude and angle, which, together with the network impedance, produces the load flows that are known to be correct at the system terminals. Power flow study is a steady state analysis used to solve such unknown values in order to plan ahead and account for various hypothetical situations of networks \cite{}.

These power flow equations have been posing much difficulty to the engineers due to their inherent nonlinearity and high computational cost of finding solutions with conventional approaches \cite{HungTuritsyn}. Iterative methods such as Newton-Raphson or Gauss-Seidel methods are the most widely used by power flow solvers for the following reasons. First, the solution can be found after several iterations in most of the cases, and the iterative methods can be easily scaled to bulk systems. Second, the solutions obtained can be used for control purposes. However, these iterative methods may be subject to divergence due to either bad initial guesses, proximity to solvable boundaries, or the nonexistence of solutions \cite{JEET}.

%Two traditional iterative algorithms for solving these non-linear power flow equations are the Newton-Raphson method and Fast Decoupled methods; however, despite their popular usage, the former method requires an expensive LU factorization of power flow Jacobian at each Newton-iteration, and the latter may subject to divergence because of either bad initial guess, proximity to solvable boundaries, or the nonexistence of solutions.
%unless the power system at hand is well-behaved. 
In order to overcome such handicaps, numerous attempts have been made to design more efficient and accurate ways to solve or certify existence of solutions to those equations \cite{coffrin2014linear,bergen, Cvijic2012}. Among them, in \cite{bolognani2014existence}, Bolognani proposed a method of certifying the solvability of power flow equations using the Banach Fixed Point Theorem. The resulting criterion has a simple algebraic representation and dramatically reduces the computational of assessing system security. The downside of the approach is the conservativeness of the criterion and its inability to certify the solvability of the whole region. The motivation for this paper is to extend Bolognani’s criterion and alleviate the conservativeness observed in the original definition.

The key contributions of this paper are as follows. We analyze an extension of Bolognani's criterion and develop a new mathematical criterion for solvability certification. Then, we derive a simple algebraic representation and propose a new technique for fast and non-iterative characterization and visualization of the solvability boundaries in a multidimensional parameter space. Finally, two important applications for solvable boundary approximation and security constraint optimization are presented.

%such that the solution boundary may be manipulated as desired;
%; one, revisit Bolognani’s method of establishing a practical method for nonlinear power flow equations and propose a modification of its criterion such that the solution boundary may be manipulated as desired; two, characterize the result from this proposed modification made.
In section \ref{sec:model}, we describe the nonlinear power flow equations and the notations that be used throughout this paper, followed by a brief explanation of the criterion established by Bolognani which is the motivation for this paper. The method for extending Bolognani’s criterion, along with its proof, be presented in section \ref{sec:modify}. Then, section \ref{sec:simulation} displays the $2$-bus simulation results and related empirical observation that discuss the relationship between the resulting ellipsoids and rhombuses. This section also discusses the algorithm, or the convex hull method for constructing outermost rhombuses that can be observed from simulations. Section \ref{sec:outermost} provides with more simulations of $3$-bus and $13$-bus radial systems, which is followed by section \ref{sec:discussion} where main results are discussed. In section \ref{sec:apps}, we introduce two possible application of this modification.

\section{Power flow problem and solvable boundary approximation} \label{sec:model}

\subsection{Power Flow Equations}

In this work we use the notation introduced originally by Bolognani in \cite{bolognani2014existence}. %The considered networks are radial ones. 
Let $\mathcal{L}:=\{1,...,n\}$ be the index set denoting all load buses, and let index $0$ correspond to the slack bus. For $h\in\mathcal{L}$, let $v_h\in\mathbb{C}$ and $i_h\in\mathbb{C}$ be the voltage and current, respectively, that define each load bus, and let $v,i\in\mathbb{C}^{n+1}$ be complex vectors with entries of $v_h$ and $i_h$, respectively. Also, let $s_h\in\mathbb{C}$ be the imposed power at each load bus with $s\in\mathbb{C}^{n+1}$ as the vector with $s_h$ as entries. We assume that the voltage level at slack bus is fixed at $v_0 = V_0$.

%For the sake of simplicity, assume that the system is in steady state and that all voltages and currents are sinusoidal signals at the same frequencies. Then each signal can be characterized as a complex number of the form $y=|y|\exp(j\angle y)$; hence 
% For the slack bus, we have that:

% \begin{equation*}
% v_0=V_0\exp(j\theta_0)
% \end{equation*}
% where $V_0, \theta_0 \in \mathbb{R}$; normally, $V_0=1$ and $\theta_0=0$.

% Then the objective of power flow equations is to solve for the voltages and currents as functions of the slack bus voltage and power at each bus. In short,
% \begin{align*}
% v_h&=v_h(V_0,\theta_0,s_1,\dots,s_n)\\
% i_h&=i_h(V_0,\theta_0,s_1,\dots,s_n)
% \end{align*}
% The construction is described below.

The power flow equations establish the relation between the power flows injected and consumed at each bus and the corresponding voltage and current levels:
%\indent For $PQ$ buses, the following holds:
\begin{equation*}
v_h\bar{i}_h=s_h\qquad \forall{h}\in\mathcal{L}:={1,\dots,n}
\end{equation*}
where the bar denotes the complex conjugate of the given vector. The power flow equations need to be supplemented with Kirchhoff laws that can be represented in the following form:
\begin{equation*}
v_L=v_0\mathbb{1}+Zi_\mathcal{L}\text{, where } Z:=Y_{\mathcal{L}\mathcal{L}}^{-1}
\end{equation*}
where $\mathbb{1}$ is a vector with all components equal to $1$, i.e. $\mathbb{1}_h = 1$.

In this work, we only consider $PQ$ buses, so the values of $s$'s are assumed to be given. The extension for $PV$ buses is rather straightforward and is also described in \cite{bolognani2014existence}.

%Following Bolognani's setup, the power flow equations are reduced to the following:

\subsection{Bolognani's certificate of solvability}

A sufficient condition for existence of a particular solution to the nonlinear power flow equations was derived in \cite{bolognani2014existence} and forms the foundation of the algorithm proposed in this work. The criterion is based on reformulation of power flow equation as fixed point of nonlinear map, i.e. equation of the form $f = G(f)$, where the map is defined via the following function:
\begin{equation}\label{map}
G(f):=-\frac{1}{V^2_0}\text{diag}(f+s_\mathcal{L})Z(\bar{f}+\bar{s}_{\mathcal{L}})
\end{equation}
and the quantity $f$ is related to bus voltage and current levels via the following relation:
\begin{equation*}
f:=v_0\bar{i}_{\mathcal{L}}-s_{\mathcal{L}}=V_0\bar{i}_{\mathcal{L}}-s_\mathcal{L}
\end{equation*}
% {\color{red} What is the definition of f ? Please replace the previous equation with the definition of f in terms of voltage and currents}

The existence of solution is certified by proving that the function $G(f)$ is contractive, and thus the solution exists. The final criterion can be compactly expressed as the following condition on the power flow levels:
\begin{equation}\label{saverio}
4\|Z\|^*_2\cdot\|s_\mathcal{L}\|_2 \leq V_0^2,
\end{equation}
with the ``nuclear'' matrix norm defined as 
\begin{align}\label{two_norm}
 \|A\|_2^*& = \max_h\sqrt{\sum_j{A_{hj}^2}}
\end{align}
Similar relation can be derived also for $1$-norms. In this case one has 
\begin{align}\label{1normcrit}
4\|Z\|^*_\infty\cdot\|s_\mathcal{L}\|_1 &\leq V_0^2 \\
 \|A\|_\infty^*& = \max_{h,k}|A_{hk}|\label{one_norm}
\end{align}
Whenever the relation \eqref{saverio} is satisfied for some power flow injection vector $s$, the power flow equations a proven to have a solution. Some additional properties of this solution were proven in the original work \cite{bolognani2014existence}.

\section{Extension of the solvability certificate} \label{sec:modify}
The solvability certificate defined by \eqref{saverio} is not unique. In fact, it was noted in the original work \cite{bolognani2014existence}, that it can be generalized by rescaling the power and voltage vectors in the original power flow equations. In this work we explore this freedom and show that these rescalings can be used to extend the criterion, and make it both less conservative and more computationally tractable.

To extend the criterion we replace $f$ with $\Lambda\varphi$, where $\Lambda$ is some diagonal matrix in $\mathbb{R}^{n \times n}$ and $\varphi$ some vertical vector in $\mathbb{R}^n$. Then we have:

\begin{equation*}
f =\Lambda\varphi =  \begin{bmatrix}\lambda_1\varphi_1\\ \lambda_2\varphi_2\\ \vdots\\ \lambda_n\varphi_n\end{bmatrix}
\end{equation*}
It follows from \eqref{map} that
\begin{equation}
\varphi_i=-\frac{1}{V_0^2}(\varphi_i+\lambda_i^{-1}s_i)\sum\limits_{j=1}^n\left(Z_{ij}\lambda_j\right)(\bar{\varphi}_j+\lambda_j^{-1}\bar{s}_j)
\end{equation}
Notice that this equation has the same structure as \eqref{map}, so the same techniques that were used to prove the solvability based on \eqref{saverio} can be applied to this equation. Then the revised criterion can be rewritten as 
\begin{equation} \label{rescaled}
4\|Z\Lambda\|^*\cdot \|\Lambda^{-1} s_\mathcal{L}\| \leq V_0^2
\end{equation}
The key idea behind the approach explored in this work is to characterize the whole union of certificates defined by the inequalities \eqref{rescaled} for all possible diagonal matrices $\Lambda$. More precisely we are interested in characterization of the set defined formally as 
\begin{align}\label{pset}
 \mathcal{P} = \left\{s_\mathcal{L}: \,\, \exists \Lambda\textrm{ such that } 4\|Z\Lambda\|^*\cdot \|\Lambda^{-1} s_\mathcal{L}\| \leq V_0^2\right\}
\end{align}
In the following sections we present the results of our analysis. In the section \ref{sec:simulation} we report the results of numerical analysis on various models, whereas in  section \ref{sec:outermost} we provide preliminary results from theoretical characterization of the set $\mathcal{P}$.

\section{Simulation results} \label{sec:simulation}
We start our analysis by considering simplest possible radial system composed of $2$ and $3$ buses. To visually characterize the union of certificates generated by different matrices $\Lambda$ in \eqref{rescaled} we can simply plot the boundaries corresponding to different $\Lambda$ on top of each other.

\begin{figure}[t!]
\centering
\includegraphics[width=0.9 \columnwidth]{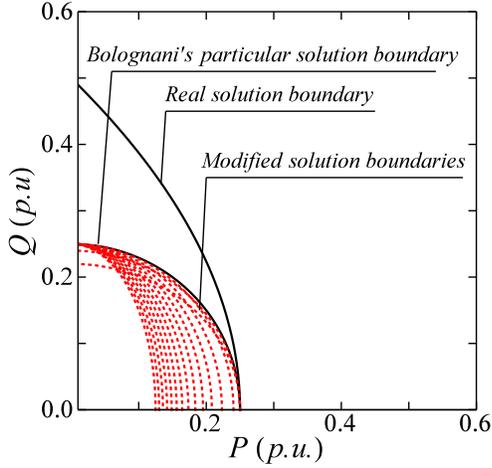}
\caption{The dotted curves illustrate the solvable boundary from varying rescalings for a $2$-bus network. The rescaled regions cover only fractions of the solution region certified by Bolognani's criterion, hence represents more conservative solution boundaries.}
\label{twobus}
\end{figure}

In Figure \ref{twobus} we present results of analysis of a simple two bus system with one slack bus and one $PQ$ bus characterized by levels of power consumption defined by $P$ and $Q$. In this model we assume, without loss of generality, that the line is purely resistive. More realistic situations can be recovered by rotating the vector $s$ in complex plane by an angle defined by the impedance of the lines. We overlay the plots of the actual solution boundary and the certificates constructed with the original criterion \eqref{saverio} as well as the modified ones \eqref{rescaled}. Several observations can be made from Figure \ref{twobus}. First, it can be clearly seen that certificates defined by \eqref{saverio} and \eqref{rescaled} do not cover the whole region where the solution to power flow equation exists. Second, one can see that rescaling does affect the region certified by the criterion. Finally, one can see that the original criterion \eqref{saverio} is the least conservative boundary, so rescaling does not increase the size of the region certified to have a power flow solution. As we will see, this last observation does not apply to more general systems discussed below.

\begin{figure}[t!]
\centering
\includegraphics[width=0.9 \columnwidth]{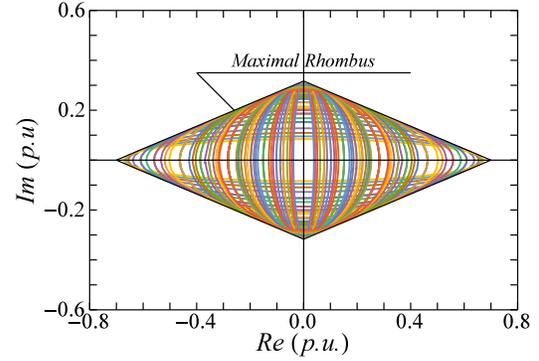}
\caption{64 elliptical solvable regions for the given 3-bus network. Note that the contour of the union of the ellipses forms a geometry identical to that of the union of rhombuses in Figure \ref{64rhombus}.}
\label{64ellipse}
\end{figure}

\begin{figure}[t!]
\centering
\includegraphics[width=0.9 \columnwidth]{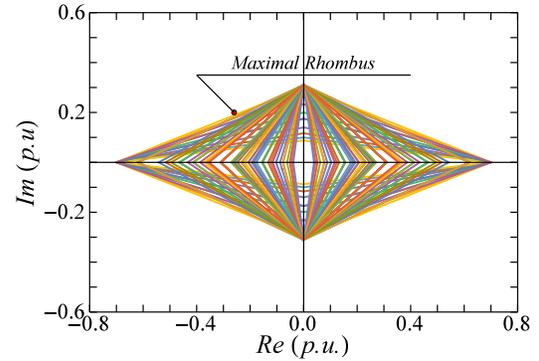}
\caption{64 rhombus-shaped solvable regions for the 3-bus network.}
\label{64rhombus}
\end{figure}

Next, we consider a more complicated case of a $3$-bus radial system with two $PQ$ buses. In the $3$-bus system, both of the power lines between buses $1-2$ and $2-3$ have the same values of resistance and inductance: $r=0.0734\, p.u.$ and $x=0.2581\, p.u.$. We consider the case in which $P_1=P_2$ and $Q_1=Q_2$. To understand the impact of the rescaling we consider diagonal matrices $\Lambda$ with elements in the range $\Lambda_{ii} \in [0.5, 25]$. The impedance matrix $Z$ in this problem is given by
\begin{equation*}
Z=\begin{bmatrix}-\frac{1}{7}&-\frac{1}{7}\\-\frac{1}{7}&-\frac{64}{203} + j\frac{2}{29}\end{bmatrix}
\end{equation*}
Figures \ref{64ellipse} and \ref{64rhombus} illustrate unions of $P$ \textit{v.s.} $Q$ plots with regions corresponding to certificates produced by different diagonal matrices $\Lambda$. Figure \ref{64ellipse} represents the union of 64-many $PQ$ plots corresponding to $||\cdot||_2^*$ norm as  and Figure \ref{64rhombus} represents the union of plots from $||\cdot||_\infty^*$ norm. The key observation that one can draw from these plots is that the outer boundary of the union of all certificates has a simple polytope, that we refer to as a ``Rhombus'' in our figures due to the shape of its projection. Second, as one can see the outer boundary improved the certificate provided by each individual $\Lambda$ and provides a less conservative way of checking the system solvability. Finally, it is important to note, that the outer boundary of both the $||\cdot||_\infty^*$-norm union and $||\cdot||_2^*$-norm union is the same. This observation suggests that this boundary can be characterized by a simple algebraic expression. 

\begin{figure}[t!]
   \centering
  \subfigure{\includegraphics[width=4.3cm]{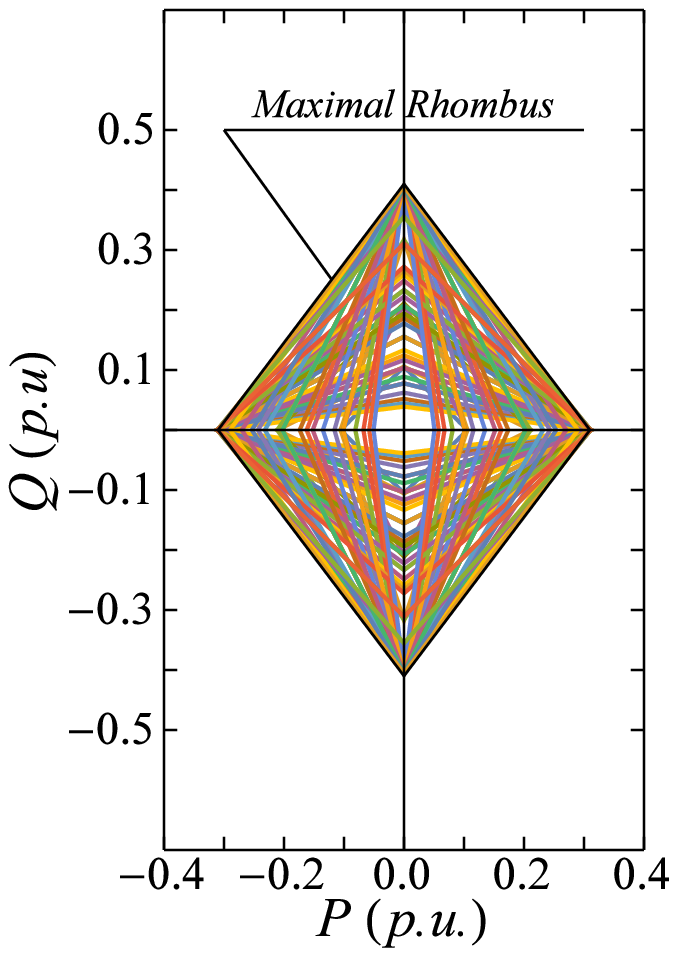}}
 \subfigure{\includegraphics[width=4.3cm]{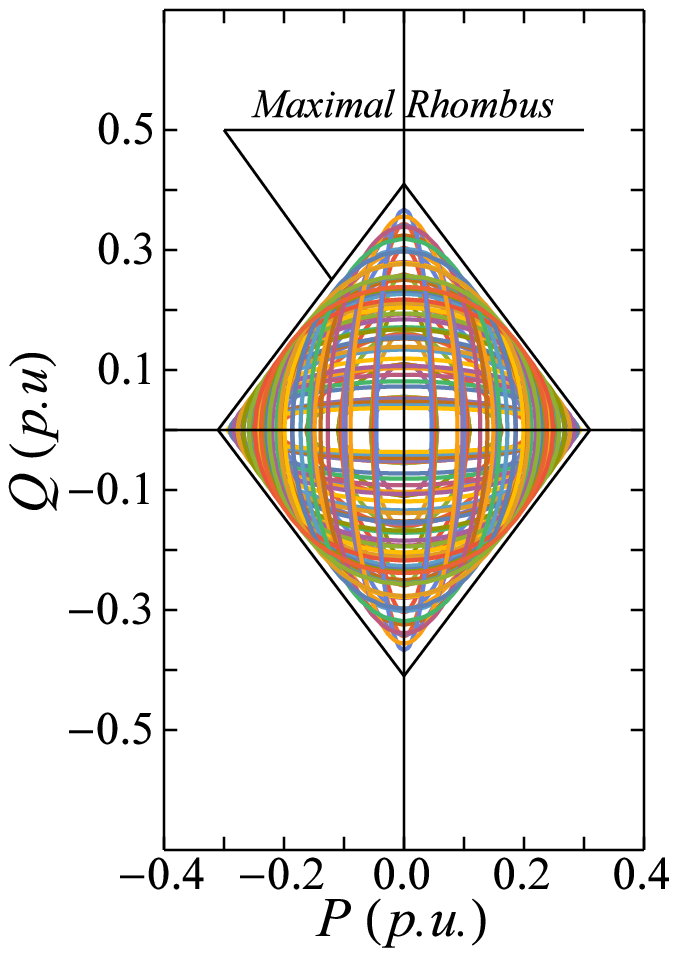}}
  \caption{Projection of $1^{st}$ and $2^{nd}$ variables' plane for 13-bus network. On the left is the rescaled regions from the definition from (\ref{one_norm}), and on the right is from (\ref{two_norm}).}%
  \label{1and2}
\end{figure}

%Figure \ref{3and10},
On figure  \ref{1and2} we provide yet another empirical confirmation of our observation by looking at the projections plots of the 13-bus network, assuming $P_1=P_2=P_3=P_4=P_5=\frac{1}{2}P_6=\frac{1}{2}P_7=\frac{1}{2}P_8=\frac{1}{2}P_9=\frac{1}{2}P_{10}=\frac{1}{4}P_{11}=\frac{1}{4}P_{12}$ and $Q_1=Q_2=Q_3=Q_4=Q_5=\frac{1}{2}Q_6=\frac{1}{2}Q_7=\frac{1}{2}Q_8=\frac{1}{2}Q_9=\frac{1}{2}Q_{10}=\frac{1}{4}Q_{11}=\frac{1}{4}Q_{12}$. The $13$-bus network configuration and data is described in \cite{nguyen2014voltage}. Moreover, fixed shunt capacitors are installed at buses $3$, $4$, $6$, $9$, $10$, $13$ with the capacities of $0.2$, $0.1$, $0.05$, $0.05$, $0.05$, $0.1$, respectively. The unit of the capacities are in $p.u.$. The impedance matrix $Z\in\mathbb{R}^{12\times12}$ is then computed, and the plot for $P$ \textit{v.s.} $Q$ is plotted using the modified criterion, using both norms $||\cdot||^*_2$ and $||\cdot||^*_\infty$. Again, we see that the union is always a simple geometric shape that looks like a rhombus on the projection plane. %{\color{red} Is it really true ? Please check convergence in figure 5. Or let's remove it in this paper and return to the question later. I'm also not sure why the maximal rhombus is exactly the same on Fig. 4 and 5. There is no reason why the polytope should have the same dimension in these two directions.} %Notice that in this particular case, the projection of $1^{st}$ and $2^{nd}$ variables' plane yields the largest outermost rhombus, hence in any dimension, it is sufficient to find the largest outermost rhombus of any projection.

%, the base power and base voltage are $S_{base} = 100\,MVA$, $V_{base} = 4.16\, kV$.

% \begin{figure}[!ht]
%     \centering
%     \subfigure{{\includegraphics[width=5cm]{r3.eps}}}%
%     \subfigure{{\includegraphics[width=5cm]{e3.eps} }}%
%     \caption{Projection of $6^{th}$ and $12^{th}$ variables' plane for 13-bus network}%
%     \label{6and12}
% \end{figure}

In the following section, based on our empirical observations we introduce a conjecture that defines the solvability certificate in terms of of the simple piecewise linear function defining the convex hull of the critical points. 

\section{Constructing the outermost rhombus} \label{sec:outermost}

In this section, we propose a method for constructing an algebraic representation of the ``Maximal rhombus'' shape that we observed in our numerical experiments. This shape is most easily constructed using the representation \eqref{1normcrit}. Consider the rescaling defined by $\lambda_k^{-1} = \max_h |Z_{hk}|$. In this case the nuclear norm of the matrix $\|Z\Lambda\|_\infty^* = 1$ that can be easily checked from the definition \eqref{one_norm}. The certificate is then reduced to 
\begin{equation} \label{convexhull}
 \sum_k \left|\frac{s_k}{s_k^{\max}}\right| \leq 1
\end{equation}
where the value of $s_k^{\max}$ is given by the following expression:
\begin{equation}
 s_k^{\max}  = \frac{V_0^2}{4 \max_h |Z_{hk}|}
\end{equation}
Note that the expression \eqref{convexhull} defines the convex hull of critical points $\mathbf{s}^{(k)} = s_k^{\max} \mathbf{e}_k$, where the unit vectors $\mathbf{e}_k$ define the rectangular coordinate system in the $s$-space. One can easily show that the polytope defined by \eqref{convexhull} is the envelope of union of polytopes defined by the rescaled version of \eqref{1normcrit}. It is also easy to check that it encompasses all the critical ellipsoids defined by the rescaled version of \eqref{saverio} with $\Lambda_k \to \infty$. At this stage we don't know the conditions under which the union of all the rescaled ellipsoids coincides with the convex hull defined by \eqref{convexhull}, however our numerical experiments show that it is often the case.

\section{Comparison with real boundary} \label{sec:discussion}

% Figure \ref{convexhull} illustrates the procedure described above for the $3$-bus case. Points A and B indicate the two endpoints of the extreme ellipse as $\lambda_2 \rightarrow 0$, and points C and D indicate the two endpoints of the extreme ellipse as $\lambda_1 \rightarrow 0$. Taking the convex hull of point A, B, C, D produces the rhombus identical to the outermost rhombus from Figure \ref{64rhombus}.

% \begin{figure}[!ht]
% \centering
% \includegraphics[width=0.7 \columnwidth]{convexhull.eps}
% \caption{Convex hull of the four endpoints}
% \label{convexhull}
% \end{figure}

\begin{figure}[t!]
\centering
\includegraphics[width=0.9 \columnwidth]{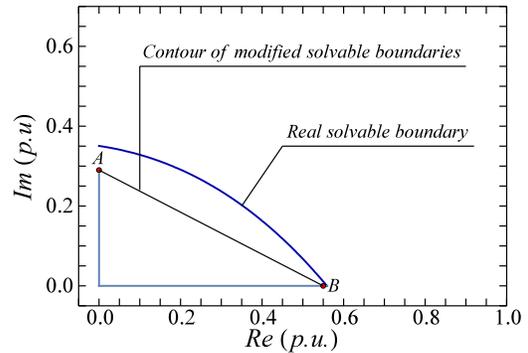}
\caption{The plot displays two different solvable boundaries of $3$-bus system; one, the actual solution boundary for a 3-bus network, and two, the contour of the union of elliptical rescaled solvable boundaries. The contour was plotted using the method for constructing the "Outermost Rhombus".}
\label{realboundary}
\end{figure}

As the certificate can be naturally used to guarantee the existence of power flow equation solutions it is reasonable to assess its conservativeness by comparing the certified region to the actual solution boundary. On Figure \ref{realboundary} we perform this comparison for the 3-bus example system. Here the outermost curve indicates the actual solution boundary found via Groebner basis approach explained in \cite{HungTuritsyn}. The polytope shown represents the first quadrant part of the convex hull as defined in \eqref{convexhull}. As expected, the polytope lies inside the solvability region, and provides a remarkably tight approximation of this boundary. %{\color{red} I suggest to remove the references to Vertical and Horizontal ellipse references from the figure, as we don't use them anymore in our construction }{\color{green} From Essie: I got rid of the two lines (extreme ellipses). Let me know if this was not the correct modification.}% Lastly, the proposed convexhull method was implemented in order to find the extreme points in the first quadrant. One can see that not only is this new linear boundary contained inside the actual solution boundary, but also that a new solution boundary can be found around a desired solution. 

% {\color{red} Hung, insert the right citation please, also correct the bib file to refer to our conference paper}

\begin{figure}[ht]
\centering
\includegraphics[width=0.85
\columnwidth]{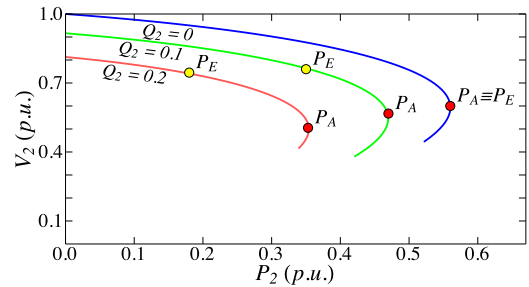}
\caption{The $PV$ curves at load bus $2$ with different levels of reactive power in the $3$-bus system.}
\label{PVcurve}
\end{figure}

As a supplement to Figure \ref{realboundary}, Figure \ref{PVcurve} presents $3$ $PV$ curves corresponding to different fixed values of reactive power $Q_2=Q_3$. $P_A$ and $P_E$ indicate the actual value and estimated value of maximum active power, respectively. For $Q_2=0$, the two values match as shown at point $B$ in Figure \ref{realboundary}. Interestingly enough, the certificate is tight for purely active power flows or DC power flow, but is more conservative as the amount of reactive power consumption is increased. %As one can observe, at $P_A$ the system encounters saddle-node bifurcation.

\section{Practical applications} \label{sec:apps}
In this section, we propose two important applications of the solvability certificate derived above.

\subsection{Non-iterative approximation of solvable boundaries}
In steady state analysis or power flow problem, the solvable boundary consists points corresponding to saddle-node bifurcation where the Jacobian matrix becomes singular. Knowing the solvable boundary of the system is essential to system operators because it indicates the limit of the system and the threats to viability of the system operation \cite{overbye1994power}. However, constructing the solvable boundary in multi-dimensional parameter space is difficult due to the complexity of the system \cite{venkatasubramanian1995local}. It has been drawing much attention for decades, and most of the proposed techniques are computationally prohibitive iterative methods \cite{overbye1994power,overbye1995computation,makarov2000computation, Chiang2008continuation}, which may suffer from the divergence problem. On the other hand, the simple linear criterion derive in this work a reasonable alternative to characterize non-iteratively the solvable boundary. Besides computational efficiency provided, a prompt technique for visualization of the solvable boundary is extremely helpful not only in contingency analysis and planning processes, but also in emergency controls where faster decision-making is required. Moreover, by construction, the approximated region bounded by the outermost rhombus is  inscribed in the actual solvability boundary; hence, the approximated boundary is closer to the feasible region in which the operational constraints are satisfied. %We envision that the approximated region is the intermediate one between the actual and the feasible ones. In this sense, the approximated boundary is more meaning to the system operators than the actual solvable boundary.

%Despite of unavoidable errors due to approximation which, however, in our test cases, are rather small;  

\subsection{Solvability criterion for security constraints}
Another important application can be developed from simple the solvability criterion described in \eqref{convexhull}. The simple criterion can be incorporated as a security constraint in planning and operational problems such as OPF. Whenever the solvability criterion is satisfied, the system is guaranteed to possess a solution to the corresponding power flow problem. The simple criterion \eqref{convexhull} is also useful for corrective control. For example, one can use this to find optimal load shedding.

Moreover, note that the proposed simple criterion can be expressed in terms of the slack bus voltage, $V_0$. Other forms of solvable criteria are described in \cite{thorp1986reactive, molzahn2012sufficient}. In the future work, we plan to extend the proposed criterion to take in account the feasibility of the solutions, including voltage and current criteria.  We expect that the new criterion could be applied directly to assess quasi-static voltage stability and would help to avoid voltage collapse in static scenarios.

\section{Conclusions}

In this work, we extended the solvability criterion derived in \cite{bolognani2014existence} by exploring the freedom associated with rescaling matrices $\Lambda$. We observed, that the union of the regions described by these certificate has often a simple polynomial shape and can be represented as a convex hull of relatively small number of points represented by the expression \eqref{convexhull}. This simple solvability certificate  depends only on the system configuration and can be naturally used as voltage security constraint.

In future works, we plan to develop a new method to identify feasible load change and design secondary voltage control based on the definition of the polytopic solvability certificate. Moreover, we plan to explore the possibility of using this criterion as static voltage stability indicator.

% \section{Meeting from 10/14/2014 Paper Outline}

% \noindent Introduction

% \begin{itemize}
%     \renewcommand{\labelitemi}{$\cdot$}
% \item Penetration of renewables
% \item Voltage stability problem in presence of renewables
% \item Previous studies: Overbye, Conejo, Ilic, Bolognani
% \item Our contribution: computationally tractable extended solvability certificates
% \end{itemize}

% \noindent Model
% \begin{itemize}
%     \renewcommand{\labelitemi}{$\cdot$}
% \item Power flow equations + motivations
% \item Bolognani's bound: review
% \item Extension ($\Lambda$) - set of certificates
% \item Union over many $\Lambda$'s $\rightarrow$ envelopes
% \end{itemize}

% \noindent Simulations for 2-bus system
% \begin{itemize}
%     \renewcommand{\labelitemi}{$\cdot$}
% \item Ellipsoids, Rhombus 
% \end{itemize}
% \noindent Algorithm for constructing outermost Rhombus
% \begin{itemize}
%     \renewcommand{\labelitemi}{$\cdot$}
% \item Derivation
% \item Prove that every ellipsoid is inside convex hull
% \end{itemize}
% \noindent Simulations
% \noindent Conclude
% \begin{itemize}
%     \renewcommand{\labelitemi}{$\cdot$}
% \item New technique for voltage stability assessment 
% \end{itemize}

% \section{Power flow approximation}
% The summation of Bolognani's paper goes here.

% \section{Solvability boundary approximation}

% Put your $\Lambda$ modified criterion there.

% \section{Simulations}

% \subsection{Solvability boundary approximation}
% \subsection{Feasible load change and secondary voltage control}
% \subsection{Probability of losing solvability}

% \section{Discussions}

\section{Acknowledgement}

We thank the NSF, MIT/Skoltech, and Masdar Initiative, as well as the Vietnam Education Foundation for their supports.

\bibliographystyle{IEEEtran}
\bibliography{main.bbl}
\end{document}